\input amstex
\magnification=\magstep1
\baselineskip=13pt
\documentstyle{amsppt}
\vsize=8.7truein
\NoRunningHeads
\def\wi{\operatorname{width}}

\def\inte{\operatorname{int}}
\def\xx{{\bold x}}
\def\yy{{\bold y}}
\def\zz{{\bold z}}
\def\ee{{\bold e}}
\def\td{\operatorname{td}}

\def\im{\operatorname{im}}
\topmatter
\title Short Rational Generating Functions for
Lattice Point Problems \endtitle
\author Alexander Barvinok and Kevin Woods \endauthor
\address Department of Mathematics, University of Michigan, Ann Arbor,
MI 48109-1109 \endaddress
\email barvinok$\@$umich.edu, kmwoods$\@$umich.edu  \endemail
\date November 2002 \enddate
\thanks This research was partially supported by NSF Grant DMS 9734138.
\endthanks
\abstract We prove that for any fixed $d$ the generating 
function of the projection of the set of integer points in a rational 
$d$-dimensional polytope can be computed in polynomial time.
As a corollary, we deduce 
that various interesting sets of lattice points,
notably integer semigroups and (minimal) Hilbert bases of rational cones, have 
short rational generating functions provided certain parameters 
(the dimension and the number of generators) are fixed. It follows then 
that many computational problems for such sets (for example, finding 
the number of positive integers not representable as a non-negative 
integer combination of given coprime positive integers $a_1, \ldots, a_d$) 
admit polynomial time algorithms. We also discuss a related problem of 
computing the Hilbert series of a ring generated by monomials. 
\endabstract
\endtopmatter
\document

\head 1. Introduction and Main Results \endhead

Our main motivation is the following question which goes back to Frobenius 
and Sylvester.
\subhead (1.1) The Frobenius Problem \endsubhead
Let $a_1, \ldots, a_d$ be positive coprime integers and let 
$$S=\Bigl\{ \mu_1 a_1 + \ldots + \mu_d a_d: \quad \mu_1, \ldots, \mu_d 
\in {\Bbb Z}_+ \Bigr\}$$
be the set of all non-negative integer combinations of $a_1, \ldots, a_d$,
or, in other words, the semigroup $S \subset {\Bbb Z}_+$ of non-negative 
integers generated by $a_1, \ldots, a_d$. 
What does $S$ look like? 
In particular, what is the largest integer not in $S$? 
(It is well known and easy to see that all sufficiently large 
integers are in $S$). How many positive integers are not in $S$? 
How many positive integers within a particular interval or a particular 
arithmetic progression are not in $S$?
\bigskip
One of the results of our paper is that for any fixed $d$ ``many'' of 
these and similar questions have ``easy'' solutions. 
For some of these questions,
notably, how to find the largest integer not in $S$, an efficient solution 
is already known \cite{K92}. 
For others, for example, how to find the number of 
positive integers not in $S$, an efficient solution was not previously known.

With a subset $S \subset {\Bbb Z}_+$ we associate the
 {\it generating function} 
$$f(S; x)=\sum_{m \in S} x^m.$$
Clearly, the series converges for all $x$ such that $|x|<1$.
We are interested in finding a ``simple'' formula for $f(S; x)$.
\example{(1.2) Examples: $d=2$ and $d=3$} Suppose that $d=2$, that is, $S$ is 
generated by two coprime positive integers $a_1$ and $a_2$. It is not 
hard to show that 
$$f(S; x)={1 -x^{a_1 a_2} \over (1-x^{a_1})(1-x^{a_2})}.$$

Suppose that $d=3$, that is, $S$ is generated
by three coprime positive integers $a_1$, $a_2$ and $a_3$. Then there 
exist (not necessarily distinct) non-negative integers 
$p_1, p_2, p_3, p_4$ and $p_5$, which can 
be computed efficiently from $a_1, a_2$ and $a_3$, such that 
$$f(S; x)={1-x^{p_1}-x^{p_2}-x^{p_3}+x^{p_4}+x^{p_5} \over 
(1-x^{a_1})(1-x^{a_2})(1-x^{a_3})}.$$
This interesting fact is, apparently, due to G. Denham \cite{D96}. 
For example, if $a=23$, $b=29$ and $c=44$, then (thanks to a MAPLE program 
written by J. Stembridge), $p_1=161$, $p_2=203$, $p_3=220$, $p_4=249$ and 
$p_5=335$. 

The idea of Denham's proof is to interpret $f(S; x)$ as the Hilbert series 
of a graded ring $M={\Bbb C}[t^{a_1}, t^{a_2}, t^{a_3}]$.
This ring $M$ can be considered as a graded module over the 
polynomial ring $R={\Bbb C}[x_1, x_2, x_3]$ 
graded by $\deg x_i=a_i$ and acting on $M$ by $x_i t^{a_j}=t^{a_i+a_j}$. Since 
the projective dimension of $M$ is 2, the Hilbert-Burch Theorem allows 
us to construct explicitly a projective resolution of $M$ and then to 
compute the Hilbert series from it, cf. Section 20.4 of \cite{E95}.

We also note that a slightly weaker form of this result is obtained by 
elementary methods in \cite{SW86}.
\endexample
What happens for $d=4$ (or larger)? Clearly, since $S$ contains all 
sufficiently large numbers, $f(S; x)$ is a rational function of the type 
$$f(S; x)=p_N(x)+{x^{N+1} \over 1-x}, \tag1.3$$
where $N$ is the largest 
integer not in $S$ and $p_N(x)$ is a polynomial of degree $N$. 
Can we find a shorter formula for $f(S; x)$?

We need some standard definitions from 
computational complexity theory (see, for example, \cite{P94}).
\definition{(1.4) Definitions} 
We define the {\it input size} of an integer $a$ as the number of 
bits needed to write $a$, that is, roughly, $1+\log_2 |a|$. Hence the 
input size of the sequence $a_1, \ldots, a_d$ will be roughly 
$d+\sum_{i=1}^d \log_2 a_i$. We are interested in the complexity of 
an algorithm which computes $f(S;x)$ from the input $a_1, \ldots, a_d$.
The algorithm is called {\it polynomial time} provided its running time 
is bounded by a certain polynomial in the input size.
\enddefinition
We show that for any fixed $d$ there is a much shorter formula for 
$f(S; x)$ than that given by (1.3).
\proclaim{(1.5) Theorem} Let us fix $d$. Then there exists a positive 
integer $s=s(d)$ and a polynomial time algorithm, which, given the 
input $a_1, \ldots, a_d$, computes $f(S;x)$ in the form
$$f(S; x)=\sum_{i \in I} \alpha_i {x^{p_i} \over (1-x^{b_{i1}}) \cdots 
(1-x^{b_{is}})},$$
where $I$ is a set of indices, $\alpha_i$ are rational numbers, 
$p_i$ and $b_{ij}$ are integers and $b_{ij} \ne 0$ for all $i,j$
\endproclaim
In particular, the number $|I|$ of fractions is bounded by a certain 
polynomial $poly$ in the input size, that is, 
in $d+\sum_{i=1}^d \log_2 a_i$.
The degree of $poly$ and the number $s=s(d)$ both grow fast with $d$,
roughly as $d^{O(d)}$. However, for any fixed $d$, the formula of 
Theorem 1.5 is much shorter than that of (1.3), in fact, 
{\it exponentially shorter}. Indeed, by \cite{EG72} it follows 
that for any fixed $d$, the integer $N$ in (1.3) can be as large as 
$O(t^2)$, where $t=\max\{a_1, \ldots, a_d\}$. Thus the length of formula 
(1.3) is quadratic in $t$, that is, {\it exponential} in the input size.
For $d=4$, there are examples (see \cite{SW86}) showing that 
if the denominator of $f(S; x)$ is chosen in the form 
$(1-x^{a_1})(1-x^{a_2})(1-x^{a_3})(1-x^{a_4})$ then the number of monomials 
in the numerator can grow as fast as $\sqrt{t}$ for 
$t=\min\{a_1, a_2, a_3, a_4\}$, which is still exponential in the input size.

Theorem 1.5 is a special case of a more general result.
Let $S \subset {\Bbb Z}^d$ be a (finite) set of integer points. For an 
integer vector $m=(\mu_1, \ldots, \mu_d)$ and (complex) variables 
$\xx=(x_1, \ldots, x_d)$, $\xx \in {\Bbb C}^d$, let 
$$\xx^m=x_1^{\mu_1} \cdots x_d^{\mu_d}$$
denote the corresponding monomial. We let $x_i^0=1$. 
Let us consider the Laurent 
polynomial 
$$f(S; \xx)=\sum_{m \in S} \xx^m.$$ 
This a priori ``long'' polynomial can sometimes be written as a ``short''
rational function 
$$f(S; \xx)=\sum_{i \in I} \alpha_i {\xx^{p_i} \over (1-\xx^{b_{i1}}) \ldots 
(1-\xx^{b_{ik}})},$$
where $\alpha_i \in {\Bbb Q}$, $p_i, b_{ij} \in {\Bbb Z}^d$ and 
$b_{ij} \ne 0$ for all $i,j$. The motivating example is the set 
$S=\bigl\{0,1,2, \ldots, n\}$, for which we have
$$f(S; x)=\sum_{k=0}^n x^n ={1-x^{n+1} \over 1-x}.$$
Thus, for this particular $S$, the long polynomial $f(S; x)$ 
can be written as a short rational function in
$x$. Indeed, writing $f(S;x)$ as a polynomial requires, roughly,
$\Omega(n \log n)$ 
bits whereas writing $f(S;x)$ as a rational function requires 
only $O(\log n)$ bits. 
A more general example is given by the set of integer points in 
a rational polyhedron.
\definition{(1.6) Definition} Let $c_1, \ldots, c_n \in {\Bbb Z}^d$ be 
integer vectors and let $\beta_1, \ldots, \beta_n \in {\Bbb Z}$ 
be integers.
The set 
$$P=\Bigl\{ x \in {\Bbb R}^d: \quad \langle c_i, x \rangle \leq \beta_i
\quad \text{for} \quad
i=1, \ldots, n \Bigr\}$$ 
is called the {\it rational polyhedron} defined by $\{c_i, \beta_i\}$. 
Again, we define the input size of $P$ as the number of bits needed to 
define $P$. That is, if $c_i=(\gamma_{i1}, \ldots, \gamma_{id})$ then 
the input size of $P$ is roughly 
$$nd+ \sum_{i=1}^n \log_2 |\beta_i| +
\sum_{i=1}^n \sum_{j=1}^d \log_2 |\gamma_{ij}|.$$
A bounded rational polyhedron is called a {\it rational polytope}. 
\enddefinition
In \cite{BP99} it is proved that for any fixed $d$, if $P \subset {\Bbb R}^d$
is a rational polyhedron 
which contains no straight lines then for $S=P \cap {\Bbb Z}^d$ the expression
$$f(S; \xx)=\sum_{m \in P \cap {\Bbb Z}^d} \xx^m$$
can be written as a short rational function. We give the precise 
statement in Theorem 3.1.  

The main result of this paper is that the 
{\it projection} of the set of integer points in 
a rational polytope has a short generating function as well.
More precisely, let $T: {\Bbb R}^d \longrightarrow {\Bbb R}^k$ be a 
linear transformation such that 
$T({\Bbb Z}^d) \subset {\Bbb Z}^k$. Thus the matrix of $T$ 
(which we also denote by $T$) with respect 
to the standard bases of ${\Bbb R}^d$ and ${\Bbb R}^k$ is integral.
The input size of $T$ is defined similarly as the number of bits needed 
to write $T$. Thus, if $T=(t_{ij})$: $i=1, \ldots, k$ and $j=1, \ldots, d$ 
then the input size of $T$ is roughly $kd +\sum_{i=1}^k \sum_{j=1}^d 
\log_2 |t_{ij}|$.
 Let $S=T(P \cap {\Bbb Z}^d)$, $S \subset {\Bbb Z}^k$, be 
the image of the set of integer points in $P$. We prove the 
following result.
\proclaim{(1.7) Theorem} Let us fix $d$. There exists a number 
$s=s(d)$ and a polynomial time algorithm, which, given a rational 
polytope $P \subset {\Bbb R}^d$ and a linear transformation 
$T: {\Bbb R}^d \longrightarrow {\Bbb R}^k$ such that 
$T({\Bbb Z}^d) \subset {\Bbb Z}^k$, computes 
the function $f(S; \xx)$ for $S=T(P \cap {\Bbb Z}^d)$,
$S \subset {\Bbb Z}^k$ in the form
$$f(S; \xx)=\sum_{i \in I} \alpha_i {\xx^{p_i} \over 
(1-\xx^{a_{i1}}) \cdots (1-\xx^{a_{is}})},$$
where $\alpha_i \in {\Bbb Q}$, $p_i$, $a_{ij} \in {\Bbb Z}^k$ and
$a_{ij} \ne 0$ for all $i,j$. 
\endproclaim
In particular, the number $|I|$ of fractions in the representation of 
$f(S; \xx)$ is bounded by a certain polynomial in the input size of 
$P$ and $T$. We do not discuss the exact dependence of 
$s(d)$ on $d$ but note that a rough estimate suggests that 
$s$ can be chosen about $d^{O(d)}$.

We obtain Theorem 1.5 as a simple corollary of Theorem 1.7 (see Section 6).
In Section 7, we discuss other interesting sets which possess short 
rational generating functions, such as the (minimal) Hilbert 
bases of rational cones and ``test sets'' 
in parametric integer programming. We also discuss a related problem 
of finding a short formula 
for the Hilbert series of a ring generated by monomials.
 
What can we do with rational generating functions? As is discussed in 
Section 3, we can efficiently perform Boolean operations on sets given
by their short rational generating functions. In particular, if 
$S_1, S_2 \subset {\Bbb Z}^d$ are two sets of integer points given
by their generating functions $f(S_1; \xx)$ and $f(S_2; \xx)$, we can 
compute the generating functions $f(S_1 \cap S_2; \xx)$, 
$f(S_1 \cup S_2; \xx)$ and $f(S_1 \setminus S_2; \xx)$ 
in polynomial time (see Theorem 3.6).
Also, by specializing at $\xx=(1, \ldots, 1)$, we can count points 
in polynomial time in finite 
sets given by their generating functions (this is not immediate 
since $\xx=(1, \ldots, 1)$ is a pole of each fraction in the 
representation of $f(S; \xx)$, cf. Theorem 2.6). 

Let $f(S; x)$ be the generating function of Theorem 1.5. Then, for the 
complement $\overline{S}={\Bbb Z}_+ \setminus S$, we compute 
the generating function $f(\overline{S}; x)=(1-x)^{-1}-f(S; x)$ 
and then compute the 
number of non-negative integers not in $S$ by specializing
$f(\overline{S}; x)$ at $x=1$. Given an interval
$[a, b] \subset {\Bbb Z}_+$, for $S'=S \cap [a,b]$,
we can compute $f(S'; x)$, and, specializing $x=1$, we can obtain the number 
of points in $S$ inside the interval $[a,b]$. 
 
The proof of Theorem 1.7 combines several methods. First, it uses 
some techniques of working with short rational generating functions 
developed by the first author, see \cite{BP99} and Sections 2 and 3. 
Second, it uses some 
``flatness''-type arguments from the geometry of numbers, see, for 
example, \cite{GLS93} and Section 4.
Finally, it relies on parametric integer programming 
arguments developed by R. Kannan, L. Lov\'asz and H. Scarf, see
\cite{K92}, [KLS90] and Section 5.
The crucial step of bringing the three ideas 
together and obtaining 
the proof of Theorem 1.7 is made by the second author (Section 6).

\remark{Remark} When a lemma or a theorem states that
``there exists a polynomial time algorithm'', the actual algorithm is 
either provided in the proof or a suitable reference is given.  
\endremark

\head 2. Rational Functions and Monomial Substitutions  \endhead

In this section, we develop certain methods of {\it specializing}
rational functions $f(\xx)$, $\xx \in {\Bbb C}^d$ 
of the type
$$f(\xx)=\sum_{i \in I} \alpha_i {\xx^{p_i} \over (1-\xx^{a_{i1}}) \cdots 
(1-\xx^{a_{ik(i)}})}, $$
where $I$ is a finite set of indices, $\alpha_i \in {\Bbb Q}$, 
$p_i, a_{ij} \in {\Bbb Z}^d$ and 
$a_{ij} \ne 0$ for all $i, j$.
We fix an upper bound $k \geq k(i)$ 
on the number of binomials in every denominator but allow
the number of variables $d$, the number $|I|$ 
of terms, the coefficients $\alpha_i$ and the vectors 
$p_i, a_{ij}$ to vary. 
Moreover, to simplify 
the notation somewhat, we will consider the case of all $k(i)$ being equal 
to a number $k$, so 
$$f(\xx)=\sum_{i \in I} \alpha_i {\xx^{p_i} \over (1-\xx^{a_{i1}}) \cdots 
(1-\xx^{a_{ik}})}. \tag2.1$$
This is a sufficiently general situation since we can always 
increase the number of binomials in a fraction by using the identity
$$\split &{\xx^p \over (1-\xx^{a_1}) \cdots (1-\xx^{a_{k-1}})}=
{\xx^p(1-\xx^{a_k}) \over (1-\xx^{a_1}) \cdots (1-\xx^{a_k})} \\=
&{\xx^p \over (1-\xx^{a_1}) \cdots (1-\xx^{a_k})} - 
{\xx^{p+a_k} \over (1-\xx^{a_1}) \cdots (1-\xx^{a_k})}. \endsplit$$
The procedure may increase the number of terms by a factor of $2^k$, but 
since $k$ is assumed to be fixed, this amounts to a constant factor increase.
 
As usual, the input size of (2.1) is the number of bits needed to 
write $f(\xx)$ down.

Let $l_1, \ldots, l_d \in {\Bbb Z}^n$ be integer vectors,
$l_i=(\lambda_{i1}, \ldots, \lambda_{in})$. The vectors define 
the {\it monomial map} $\phi: {\Bbb C}^n \longrightarrow {\Bbb C}^d$ as 
follows:
$$\aligned &\zz \longmapsto \xx \\
&(\zeta_1, \ldots, \zeta_n) \longmapsto (x_1, \ldots, x_d), \quad 
\text{where} \quad x_i=\zz^{l_i}. \endaligned \tag2.2$$
The input size of this monomial map is the number of bits needed to define 
it, that is, roughly, $dn + \sum_{i=1}^d \log_2 |\lambda_{ij}|$. 

Suppose that the image of $\phi$ does not consist entirely of poles of 
$f(\xx)$. Then we can define a rational function 
$g: {\Bbb C}^n \longrightarrow {\Bbb C}$ by 
$$g(\zz)=f\bigl(\phi(\zz)\bigr).$$
\bigskip
The goal of this section is to construct a polynomial time algorithm,
which, given a rational function (2.1) with a 
fixed number $k$ of binomials in each fraction and a 
monomial substitution (2.2),
computes a formula for $g(\zz)$.
Note that we cannot just substitute 
$\xx=\phi(\zz)$ in the formula (2.1) since any $\zz \in {\Bbb C}^n$ 
may turn out to be a pole for some fraction of (2.1) and yet a regular point 
of $g$. For example, if $d=1$, $n=0$ and 
$$f(x)={1 \over 1-x} - {x^{n+1} \over 1-x}=\sum_{m=0}^n x^m,$$
then $x=1$ is the pole of both fractions but is a regular point of $f$; 
we have $f(1)=n+1$.

To this end, let us associate with the rational function (2.1) 
a meromorphic function $F(c)$, $c \in {\Bbb C}^d$, defined by 
$$F(c)=\sum_{i \in I} \alpha_i {\exp \langle c, p_i\rangle 
\over (1-\exp \langle c, a_{i1} \rangle) \cdots 
(1-\exp \langle c, a_{ik} \rangle)}. \tag2.3$$
As usual, for $c \in {\Bbb C}^d$ with $c=r+it$, 
where $r,t \in {\Bbb R}^d$ and 
$a \in {\Bbb R}^d$, we let $\langle c, a \rangle = \langle r, a \rangle +
i \langle t, a\rangle$, where $\langle \cdot , \cdot \rangle$ is the 
standard scalar product in ${\Bbb R}^d$. The 
set of poles of the $i$-th fraction is the union
over $1 \leq j \leq k$ of the hyperplanes 
$\bigl\{c\in {\Bbb C}^d:\ \langle c, a_{ij} \rangle=0\bigr\}$.
 However, the set of 
poles of $F(c)$ may be much smaller because of cancellations of 
singularities.   

There is a simple relation 
between (2.1) and (2.3).
For $c=(\gamma_1, \ldots, \gamma_d)$ and 
$\xx=(x_1, \ldots, x_d)$ we write 
$$\xx=\ee^c \quad \text{provided} \quad x_i =\exp\{\gamma_i\} \quad 
\text{for} \quad i=1, \ldots, d.$$ 
Then the functions (2.1) and (2.3) are related by the equation
$$F(c)=f\bigl(\ee^c\bigr).$$

Let $L \subset {\Bbb  C}^d$ be a subspace such that a generic 
$c \in L$ is a regular point of $F(c)$. We want to construct a 
short formula for $F(c)$ for $c \in L$.
We assume that the subspace $L \subset {\Bbb C}^d$ is given by its integer 
basis. Again, we cannot just 
use (2.3), since $L$ may be orthogonal to some vectors $a_{ij}$ and hence
a generic $c \in L$ may be a pole of some fractions in (2.3) while 
being a regular point of $F(c)$.
\definition{(2.4) Definition} Given $l$, let us consider the function
$$G(\tau; \xi_1, \ldots, \xi_l)=\prod_{i=1}^l 
{\tau \xi_i \over 1-\exp\{-\tau \xi_i\}}$$ in $l+1$ 
(complex) variables $\tau$ and $\xi_1, \ldots, \xi_l$. It is easy to 
see that $G$ is analytic in a neighborhood of the origin
$\tau=\xi_1= \ldots =\xi_l=0$ and therefore there exists an 
expansion
$$G(\tau; \xi_1, \ldots, \xi_l)=\sum_{j=0}^{+\infty} 
\tau^j \td_j(\xi_1, \ldots, \xi_l),$$
where $\td_j(\xi_1, \ldots, \xi_l)$ is a homogeneous polynomial of degree
$j$, called the $j$-th {\it Todd polynomial} in $\xi_1, \ldots, \xi_l$.
It is easy to check that $\td_j(\xi_1, \ldots, \xi_l)$ is a
symmetric polynomial with rational coefficients, cf. \cite{BP99}.
\enddefinition
\proclaim{(2.5) Lemma} Let us fix $k$.
Then there exists a polynomial time 
algorithm, which, given a function (2.3)
and a subspace $L \subset {\Bbb C}^d$ which does not lie entirely in 
the set of poles of $F$
computes 
$F(c)$ for $c \in L$ in the form 
$$F(c)=\sum_{i \in I'} \beta_i {\exp\langle c, q_i\rangle \over 
\bigl(1-\exp \langle c, b_{i1} \rangle \bigr) \cdots 
\bigl(1-\exp \langle c, b_{is}\rangle \bigr)},$$
where $s \leq k$, $\beta_i \in {\Bbb Q}$, 
 $q_i, b_{ij} \in {\Bbb Z}^d$ 
and $b_{ij}$ is not orthogonal to $L$ for any $i,j$.
\endproclaim 
\demo{Proof} 
Let us consider the representation (2.3). Let us 
choose a vector $v \in {\Bbb R}^d$ such that 
$\langle v, a_{ij} \rangle \ne 0$ for all $a_{ij}$.
Such a vector $v$ can be constructed in polynomial time, see, 
for example, \cite{BP99}. Let $\tau$ be a complex parameter. Then, for any 
regular point $c$ of $F(c)$ the function $F(c+\tau v)$ is
an analytic function in a neighborhood of $\tau=0$ and
the constant term of its expansion at $\tau=0$ is equal to $F(c)$.
Hence our goal is to compute the constant term (in $\tau$) of every fraction 
in the representation (2.3) of $F(c +\tau v)$ and add them up.

Let us consider a typical fraction
$$h(\tau)={\exp \langle c +\tau v, p  \rangle  \over 
\bigl(1-\exp \langle c +\tau v, a_1 \rangle \bigr) \cdots 
\bigl(1-\exp \langle c + \tau v, a_k \rangle \bigr)},$$
where $p, a_j \in {\Bbb Z}^d$,
as a function of $\tau$. Suppose that 
the vectors $a_i$ orthogonal to $L$ are $a_1, \ldots, a_l$ for some 
$l \leq k$. 
Then 
$$\split h(\tau)=&\tau^{-l}\exp\langle c, p \rangle
 \exp\bigl\{\tau \langle v, p \rangle\bigr\}
\prod_{i=1}^l { \tau \over 1 -\exp\bigl\{\tau \langle v, a_i \rangle\bigr\}} \\
\times &\prod_{i=l+1}^k 
{1 \over 1-\exp \langle c+\tau v, a_i \rangle}. \endsplit$$  
Now we observe that $\tau^l h(\tau)$ is an analytic function of 
$\tau$ and that our goal is to compute the coefficient of 
$\tau^l$ in the expansion of $\tau^l h(\tau)$ in the neighborhood
of $\tau=0$.  

First, we observe that 
$$\exp\bigl\{ \tau \langle v, p \rangle\bigr\}=
\sum_{j=0}^{+\infty} {\langle v, p \rangle^j \over j!} \tau^j. \tag2.5.1$$
     
Second, letting $\xi_i=-\langle v, a_i \rangle$ for $i=1, \ldots, l$, 
we observe that 
$$\prod_{i=1}^l { \tau \over 1 -\exp\bigl\{\tau \langle v, a_i \rangle\bigr\}}
={1 \over \xi_1 \cdots \xi_l} 
\sum_{j=0}^{+\infty} \tau^j \td_j(\xi_1, \ldots, \xi_l) \tag2.5.2$$

Finally,
$$\prod_{i=l+1}^k 
{1 \over 1-\exp \langle c+\tau v, a_i \rangle}=
\sum_{j=0}^{+\infty} H_j(c, a_{l+1}, \ldots, a_k, v) \tau^j \tag2.5.3$$
for some functions $H_j$.

Note that $\tau=0$ is a regular point of 
$$\prod_{i=l+1}^k 
{1 \over 1-\exp \langle c+\tau v, a_i \rangle}$$
and so we compute $H_j$ differentiating the product $j$ times and setting 
$\tau=0$.
By the repeated application of the chain rule, 
$H_j$ 
is a polynomial in
$\exp\langle c, a_i \rangle$, $\langle v, a_i\rangle$ and 
$(1-\exp\langle c, a_i\rangle)^{-1}$. Thus,
for all $j_1, j_2, j_3$ such that $j_1+j_2+j_3=l$,
 we have to combine the $j_1$-st
term of (2.5.1) , the $j_2$-nd term of (2.5.2) and the $j_3$-rd term of 
(2.5.3).
Since $l \leq k$ and $k$ is fixed, we get the desired result.
{\hfill \hfill \hfill} \qed
\enddemo
\remark{Remark} If $L=\{0\}$ and $0$ is a regular point of $F(c)$, the 
algorithm of Lemma 2.5 computes the number $F(0)$. This procedure is used 
in \cite{B94} to compute the number of integer points in a polytope.
\endremark
Now we can compute the result of the monomial substitution (2.2) 
into the rational function (2.1).
\proclaim{(2.6) Theorem} Let us fix $k$. Then there exists a 
polynomial time algorithm, which, given a 
function (2.1) and a monomial map
$\phi: {\Bbb C}^n \longrightarrow {\Bbb C}^d$ given by (2.2), 
such that the image of $\phi$ does not lie entirely in the set of poles 
of $f(\xx)$ computes the function 
$g(\zz)=f\bigl(\phi(\zz)\bigr)$ as 
$$g(\zz)=\sum_{i\in I'} \beta_i {\zz^{q_i} \over 
(1-\zz^{b_{i1}}) \cdots (1-\zz^{b_{is}})},$$
where $s \leq k$, $\beta_i \in {\Bbb Q}$, $q_i, b_{ij} \in {\Bbb Z}^n$ 
and $b_{ij} \ne 0$ for all $i,j$. 
\endproclaim
\demo{Proof}
Let $F(c)$ be the function (2.3) associated to $f(\xx)$.
With the monomial map (2.2) we associate a 
linear transformation $\Phi: {\Bbb C}^n \longrightarrow {\Bbb C^d}$  
$$c \longmapsto \bigl( \langle c, l_1\rangle, \ldots, 
\langle c, l_d \rangle \bigr)$$
and the adjoint transformation $\Phi^{\ast}: {\Bbb C}^d \longrightarrow 
{\Bbb C}^n$,
$$\Phi^{\ast}(\xi_1, \ldots, \xi_d)=\xi_1 l_1 + \ldots + \xi_d l_d.$$

Let us define
$$G(c)=F\bigl(\Phi(c)\bigr) \quad \text{for} \quad c \in {\Bbb C}^n.$$
Hence 
$$G(c)=g(\ee^c).$$
Let $L \subset {\Bbb C}^d$ be the image of ${\Bbb C}^n$ under $\Phi$. 
Then $L$ does not lie entirely in the set of poles of $F(c)$. Applying 
Lemma 2.5, we compute $G(c)=F\bigl(\Phi(c)\bigr)$ in the form
$$G(c)=\sum_{i\in I'} \beta_i
{\exp \langle \Phi(c), u_i \rangle  \over 
(1-\exp\langle \Phi(c), v_{i1} \rangle) \cdots (1-\exp\langle \Phi(c),
v_{is} \rangle)},$$
where for $i,j$ we have $\langle \Phi(c), v_{ij} \rangle \ne 0$ for 
a generic $c \in L$.
Now we let $q_i=\Phi^{\ast}(u_i)$ and $b_{ij}=\Phi^{\ast}(v_{ij})$ 
so that 
$$g(\ee^c)=G(c)=\sum_{i\in I'} \beta_i
{\exp \langle c, q_i \rangle  \over 
(1-\exp\langle c, b_{i1} \rangle) \cdots (1-\exp\langle c,
b_{is} \rangle)}$$
and the result follows.
{\hfill \hfill \hfill} \qed
\enddemo
\remark{Remark} In particular, if $\xx=(1, \ldots, 1)$ is a regular point 
of (2.1), we can choose $l_1= \cdots =l_d=0$ in (2.2). In this case,
the algorithm of Theorem 2.6 computes the value of $f(1, \ldots, 1)$.
\endremark

\head 3. Operations with Generating Functions \endhead

Some of the results of this section are stated in \cite{BP99}. 
Many of the proofs in \cite{BP99} are 
only sketched and some non-trivial details are omitted. 
We give a mostly independent presentation with 
complete proofs. The main goal of this section is to prove that 
if finite sets $S_1, S_2 \subset {\Bbb Z}^d$ are given by their 
generating functions $f(S_1; \xx)$ and $f(S_2; \xx)$ then the generating 
function $f(S; \xx)$ of their intersection $S=S_1 \cap S_2$ can be computed
efficiently. Our main tool is the generating function for the integer 
points in a rational polyhedron.

Let $P \subset {\Bbb R}^d$ be a rational polyhedron and let 
$S=P \cap {\Bbb Z}^d$ be the set of integer points in $P$. Let 
$$f(S; \xx)=\sum_{m \in P \cap {\Bbb Z}^d} \xx^m.$$
Thus if $P$ is bounded,  $f(S; \xx)$ is a Laurent polynomial in $\xx$.
If $P$ (possibly unbounded) does not contain straight lines then 
there is a non-empty open set $U \subset {\Bbb C}^d$ such that
 the series converges absolutely and uniformly 
on compact subsets of $U$ to a rational function of $\xx$.
If $P$ contains a straight line it is convenient to agree that 
$f(S; \xx)\equiv 0$, see \cite{BP99}. 

We need the following result from \cite{BP99} which states that $f(S; \xx)$
can be written as a short rational function.
\proclaim{(3.1) Theorem} Let us fix $d$. Then there exists a polynomial 
time algorithm, which, for any given rational polyhedron
$P \subset {\Bbb R}^d$ computes $f(P\cap {\Bbb Z}^d; \xx)$ as 
$$f(P \cap {\Bbb Z}^d; \xx)=\sum_{i \in I} \epsilon_i {\xx^{p_i} \over 
(1-\xx^{a_{i1}}) \cdots (1-\xx^{a_{id}})},$$
where $\epsilon_i \in \{-1, 1\}$, $p_i, a_{ij}, \in {\Bbb Z}^d$ and 
$a_{ij} \ne 0$ for all $i, j$. In fact, for each $i$,
$a_{i1}, \ldots, a_{id}$ is a basis of ${\Bbb Z}^d$.
\endproclaim
A (complete) proof can be found in \cite{BP99}, Theorem 4.4. 

To compute the generating function of the intersection of two sets, we 
compute a more general operation, that is,
 the Hadamard product of two rational 
generating functions.
\definition{(3.2) Definition} Let $g_1$ and $g_2$ be Laurent power series in
 $\xx \in {\Bbb C}^d$
$$g_1(\xx)=\sum_{m \in {\Bbb Z}^d} \beta_{1m} \xx^m \quad 
\text{and} \quad g_2(\xx)=\sum_{m \in {\Bbb Z}^d} \beta_{2m} \xx^m.$$
The {\it Hadamard product} $g=g_1 \star g_2$ is the power series
$$g(\xx)=\sum_{m \in {\Bbb Z^d}} \beta_m \xx^m \quad 
\text{where} \quad \beta_m=\beta_{1m} \beta_{2m}.$$ 
\enddefinition
\medskip
First we will show that the Hadamard product of the Laurent expansions 
of some particular rational functions can be computed in polynomial time.
Namely, let us choose a non-zero vector $l \in {\Bbb Z}^d$ and suppose 
that $a_{11}, \ldots, a_{1k} \in {\Bbb Z^d}$ and 
$a_{21}, \ldots, a_{2k} \in {\Bbb Z}^d$ are vectors such that 
$\langle l, a_{ij} \rangle <0$ for all $i,j$.
Let $p_1, p_2 \in {\Bbb Z}^d$ and let 
$$g_1(\xx)={\xx^{p_1} \over (1-\xx^{a_{11}}) \cdots (1-\xx^{a_{1k}})}
\quad 
\text{and} \quad 
g_2(\xx)={\xx^{p_2} \over (1-\xx^{a_{21}}) \cdots (1-\xx^{a_{2k}})}.
\tag3.3$$
We observe that for all $\xx$ in a sufficiently small neighborhood $U$ 
of $\xx_0=\ee^l$, we have $|\xx^{a_{ij}}|<1$ and so $g_1$ and 
$g_2$ have Laurent series expansions for $\xx \in U$.
Indeed, if $|\xx^a|<1$, the fraction $1/(1-\xx^a)$ expands as a geometric 
series 
$${1 \over 1-\xx^a}=\sum_{\mu \in {\Bbb Z}_+} \xx^{\mu a},$$
and to obtain the expansions of $g_1$ and $g_2$ we multiply the 
corresponding series.
Clearly, the Hadamard product of the expansions converges 
for all $\xx \in U$ to some analytic function $h$, which we 
also denote $g_1 \star g_2$.
We prove that once the number $k$ of binomials in (3.3) is fixed, 
there is a polynomial time 
algorithm for computing the Laurent expansion of $h=g_1 \star g_2$ as a
short rational function.
\proclaim{(3.4) Lemma} Let us fix $k$. Then there exists 
a polynomial time algorithm, which, given functions (3.3)
such that for some $l \in {\Bbb Z}^d$ we have
$\langle a_{ij}, l \rangle <0$ for all $i,j$, computes
a function $h(\xx)$ in the form
$$h(\xx)=\sum_{i \in I} \beta_i {\xx^{q_i} \over 
(1-\xx^{b_{i1}}) \cdots (1-\xx^{b_{is}})}$$ 
with $q_i, b_{ij} \in {\Bbb Z}^d$, $\beta_i \in {\Bbb Q}$ and 
$s \leq 2k$ such that 
$h$ has the Laurent expansion in a neighborhood $U$ of 
$\xx_0=\ee^l$ 
and $h(\xx)=g_1(\xx) \star g_2(\xx)$.
\endproclaim
\demo{Proof} In the space ${\Bbb R}^{2k}=
\bigl\{(\xi_1, \ldots, \xi_{2k}) \bigr\}$ let $P$ be a rational 
polyhedron defined by the equations
$$p_1+\xi_1 a_{11} + \ldots +\xi_k a_{1k} = 
p_2+\xi_{k+1} a_{21} + \ldots + \xi_{2k} a_{2k}$$ 
and the inequalities 
$$\xi_i \geq 0 \quad \text{for} \quad i=1, \ldots, 2k.$$
Let $\zz=(z_1, \ldots, z_{2k})$ and let us consider 
the series 
$$f(P \cap {\Bbb Z}^{2k}; \zz)=\sum_{m \in P \cap {\Bbb Z}^{2k}} \zz^m.
\tag3.4.1$$
Clearly, the series converges absolutely and uniformly on compact sets
as long as $|z_i|<1$ for $i=1, \ldots, 2k$.
By Theorem 3.1 we compute $f(P \cap {\Bbb Z}^{2k}; \zz)$ in the form
$$f(P \cap {\Bbb Z}^{2k}; \zz)=\sum_{i \in I'} \epsilon_i 
{\zz^{u_i} \over (1-\zz^{v_{i1}}) \cdots (1-\zz^{v_{i(2k)}})}, \tag3.4.2$$
for some vectors $u_i, v_{ij} \in {\Bbb Z}^{2k}$ and some  
numbers $\epsilon_i\in \{-1,1\}$, where $v_{ij} \ne 0$ for all $i,j$.

On the other hand, expanding $g_1(\xx)$ and $g_2(\xx)$ as products 
of geometric series, we obtain
$$\split &g_1(\xx)=\xx^{p_1}\prod_{i=1}^k \sum_{\mu_i \in {\Bbb Z}_+} 
\xx^{\mu_i a_i}=\sum_{(\mu_1, \ldots, \mu_k) \in  {\Bbb Z}^k_+} 
\xx^{p_1 + \mu_1 a_{11} + \ldots + \mu_k a_{1k}} \quad \text{and} \\ 
&g_2(\xx)=\xx^{p_2}\prod_{i=1}^k \sum_{\nu_i \in {\Bbb Z}_+} 
\xx^{\nu_i a_i}=\sum_{(\nu_1, \ldots, \nu_k) \in {\Bbb Z}^k_+} 
\xx^{p_2 + \nu_1 a_{21} + \ldots + \nu_k a_{2k}}. \endsplit$$
Since the Hadamard product is bilinear and since
$$\xx^{m_1} \star \xx^{m_2}=\cases \xx^{m_1} &\text{if\ } m_1=m_2 \\ 0 &
\text{if\ } m_1 \ne m_2, \endcases$$
we conclude that 
$$g_1(\xx) \star g_2(\xx)=
\xx^{p_1} 
\sum \Sb (m, n) \in P \cap {\Bbb Z}^{2k} \\
m=(\mu_1, \ldots, \mu_k) \\
n=(\nu_1, \ldots, \nu_k) \endSb
\xx^{\mu_1 a_{11} + \ldots + \mu_k a_{1k}}.$$
Thus $h(\xx)$ is obtained from the function $\xx^{p_1} f(P \cap {\Bbb Z}^{2k}; \zz)$ (cf. (3.4.1)--(3.4.2)) 
by the monomial substitution 
$$z_1=\xx^{a_{11}}, \ldots, z_k=\xx^{a_{1k}}, z_{k+1}=1, \ldots, 
z_{2k}=1.$$
Now we use Theorem 2.6 to compute the result of the monomial substitution 
in (3.4.2).
{\hfill \hfill \hfill} \qed
\enddemo 

Now we are ready to prove the main result of this section.
Suppose we have two finite sets $S_1, S_2 \subset {\Bbb Z}^d$ and 
let $f(S_1; \xx)$ and $f(S_2; \xx)$ be the corresponding generating functions
$$f(S_1; \xx)=\sum_{m \in S_1} \xx^m \quad
\text{and} \quad 
f(S_2; \xx)=\sum_{m \in S_2} \xx^m.$$
Suppose further, that $f(S_1; \xx)$ and $f(S_2; \xx)$ can be written as 
short rational functions
$$\aligned 
&f(S_1; \xx)=\sum_{i \in I_1}\alpha_i {\xx^{p_i} \over (1-\xx^{a_{i1}}) 
\cdots (1-\xx^{a_{ik}})} \quad 
\text{and} \\
&f(S_2; \xx)=\sum_{i \in I_2}\beta_i {\xx^{q_i} \over (1-\xx^{b_{i1}}) 
\cdots (1-\xx^{b_{ik}})} \endaligned \tag3.5$$
with $\alpha_i, \beta_i \in {\Bbb Q}$, 
$p_i, q_i, a_{ij}, b_{ij} \in {\Bbb Z}^d$ and 
$a_{ij}, b_{ij} \ne 0$.
Now let us consider $S_1$ and $S_2$ as {\it defined} by 
representations (3.5) of 
$f(S_1; \xx)$ and $f(S_2; \xx)$ as rational functions. 
Let $S=S_1 \cap S_2$. Our goal is to compute the representation of 
$$f(S; \xx)=\sum_{m \in S} \xx^m$$ 
as a short rational function. Again, we assume the number of $k$ of binomials
in each fraction of (3.5) fixed and allow numbers $\alpha_i$ and 
$\beta_i$ and vectors $p_i,q_i$ and 
$a_{ij}, b_{ij}$ to vary.
\proclaim{(3.6) Theorem} Let us fix $k$. Then there exists a polynomial time 
algorithm, which, given $f_1(S_1; \xx)$ and $f_2(S_2; \xx)$ computes 
$f(S; \xx)$ for $S=S_1 \cap S_2$ in the form
$$f(S; \xx)=\sum_{i \in I} \gamma_i {\xx^{u_i} \over 
(1-\xx^{v_{i1}}) \cdots (1-\xx^{v_{is}})},$$
where $s \leq 2k$, $\gamma_i \in {\Bbb Q}$, $u_i, v_{ij} \in {\Bbb Z}^d$ 
and $v_{ij} \ne 0$ for all $i,j$. 
\endproclaim
\demo{Proof} Let us choose a vector $l \in {\Bbb Z}^d$, such that 
$\langle l, a_{ij} \rangle \ne 0$ and $\langle l, b_{ij} \rangle \ne 0$ 
for all $i,j$. As we remarked before, such a vector $l$ can be constructed 
in polynomial time. When $\langle l, a_{ij} \rangle >0$ 
or when $\langle l, b_{ij} \rangle >0$ we apply the identity
$${\xx^p \over 1-\xx^a}=-{\xx^{p-a} \over 1-\xx^{-a}},$$  
to reverse the direction of $a_{ij}$ or $b_{ij}$,
so that we achieve $\langle l, a_{ij} \rangle <0$ and 
$\langle l, b_{ij} \rangle <0$ for all $i,j$ in the representations
(3.5). Then we can write 
$$f(S_1; \xx)=\sum_{i \in I_1} \alpha_i g_{1i}(\xx) \quad 
\text{and} \quad f(S_2; \xx)=\sum_{i \in I_2} \beta_i g_{2i}(\xx)$$
for some functions $g_{i1}, g_{2i}$ of type (3.3). 
There are Laurent 
series expansions of $f(S_1; \xx)$ and $f(S_2; \xx)$ in a neighborhood $U$ of
the point $\xx_0=\ee^l$ and 
$$f(S; \xx)=f(S_1; \xx) \star f(S_2; \xx)=
\sum_{i_1 \in I_1, i_2 \in I_2} \alpha_{i_1} \beta_{i_2} 
g_{1i}(\xx) \star g_{2i}(\xx).$$
We use Lemma 3.4 to compute $f(S; \xx)$.
{\hfill \hfill \hfill} \qed   
\enddemo
Let $S_1, \ldots, S_m \subset {\Bbb Z}^d$ be sets. We say that 
$S \subset {\Bbb Z}^d$ is a {\it Boolean combination} 
of $S_1, \ldots, S_m$ 
provided $S$ is obtained from $S_i$ by taking intersections, unions and 
complements. An immediate corollary of Theorem 3.6 is that the generating 
function of a Boolean combination of sets can be computed in polynomial 
time.
\proclaim{(3.7) Corollary} Let us fix $m$ (the number of sets 
$S_i \subset {\Bbb Z}^d$) and $k$ (the number of binomials in each 
fraction of $f(S_i; \xx)$). Then there exists an $s=s(k,m)$ and a 
polynomial time algorithm, which, for any $m$ (finite) sets 
$S_1, \ldots, S_m \subset {\Bbb Z}^d$ given by their 
generating functions $f(S_i; \xx)$ and a set $S \subset {\Bbb Z}^d$
 defined 
as a Boolean combination of $S_1, \ldots, S_m$, computes 
$f(S; \xx)$ in the form
$$f(S; \xx)=\sum_{i \in I} \gamma_i {\xx^{u_i} \over 
(1-\xx^{v_{i1}}) \cdots (1-\xx^{v_{is}})},$$
where $\gamma_i \in {\Bbb Q}$, $u_i, v_{ij} \in {\Bbb Z}^d$ 
and $v_{ij} \ne 0$ for all $i,j$.
\endproclaim
\demo{Proof} We note that 
$$\split &f(S_1 \cup S_2; \xx)=f(S_1; \xx)+
f(S_2; \xx)-f(S_1 \cap S_2; \xx) 
\quad \text{and} \\
& f(S_1 \setminus S_2; \xx)=f(S_1; \xx)-
f(S_1 \cap S_2; \xx) \endsplit$$
for any two subsets $S_1, S_2 \subset {\Bbb Z}^d$.   
The proof follows by Theorem 3.6.
{\hfill \hfill \hfill} \qed
\enddemo
Finally, we discuss how to {\it patch} together several generating functions
into a single generating function.
\definition{(3.8) Definitions} By the {\it interior}
$\inte P$ of a polyhedron $P \subset {\Bbb R}^d$ we always 
mean the relative interior, that is, the interior of $P$ with respect to 
its affine hull.  

Let $X \subset {\Bbb R}^d$ be a set. 
We denote by $[X]$ the indicator function 
$[X]: {\Bbb R}^d \longrightarrow {\Bbb R}$, 
$$[X](x)=\cases 1 &\text{if\ } x \in X \\ 0 &\text{if \ } x \notin X. 
\endcases$$
We will need a simple formula for the indicator of the relative interior 
of a polytope:
$$[\inte P]=(-1)^{\dim P} \sum_{F} (-1)^{\dim F} [F], \tag3.8.1$$
where the sum is taken over all faces of $P$ including $P$ itself. 
This is a simple corollary of the Euler-Poincar\'e formula;
see, for example, Section VI.3 of \cite{B02}.
\enddefinition
From Theorem 3.1 we deduce the following corollary.
\proclaim{(3.9) Corollary} Let us fix $d$. Then there exists a polynomial 
time algorithm, which, for any given rational polytope 
$P \subset {\Bbb R}^d$ computes 
$f(S; \xx)$ with $S=\bigl(\inte P \bigr) \cap {\Bbb Z}^d$ in the 
form
$$f(S; \xx)=\sum_{i \in I} \alpha_i {\xx^{p_i} \over 
(1-\xx^{a_{i1}}) \cdots (1-\xx^{a_{id}})},$$
where $\alpha_i \in {\Bbb Q}$, $p_i, a_{ij} \in {\Bbb Z}^d$ and 
$a_{ij} \ne 0$ for all $i,j$.
\endproclaim
\demo{Proof} Applying formula (3.8.1), we get 
$$f(S; \xx)=(-1)^{\dim P} \sum_{F} (-1)^{\dim F} f(F \cap {\Bbb Z}^d; \xx).$$
Since the dimension $d$ is fixed, there are polynomially many faces $F$
and their descriptions can be computed in polynomial time from the description
of $P$.
We use Theorem 3.1 to complete the proof.
{\hfill \hfill \hfill} \qed
\enddemo
Let us consider the following situation. 
Let $S \subset {\Bbb Z}^d$ be a finite set and let 
$Q_1, \ldots, Q_n \subset {\Bbb R}^d$ be a
collection of rational polytopes such that
$S \subset \bigcup_{i=1}^n \inte Q_i$ and $\inte Q_i \cap \inte Q_j=\emptyset$ 
for $i \ne j$. In a typical situation, $Q_1, \ldots, Q_n$ is a polytopal 
complex, that is, the intersection of every two polytopes $Q_i$ and 
$Q_j$, if non-empty, is a common face of $Q_i$ and $Q_j$ and a face 
of a polytope $Q_i$ from the collection is also a polytope from the 
collection (in particular, not all $Q_i$ are full-dimensional). 
In this case, $\bigcup_{i=1}^n Q_i=\bigcup_{i=1}^n \inte Q_i$ and 
$\inte Q_i$ are pairwise disjoint.

 Suppose that we are given the functions
$$f(S \cap Q_j; \xx)=\sum_{i \in I_j} \alpha_{i,j} {\xx^{p_{i,j}} 
\over (1-\xx^{a_{i1, j}}) \cdots (1-\xx^{a_{ik,j}})}$$
and that we want to compute $f(S; \xx)$. In other words, we want 
to patch together several generating functions $f(S \cap Q_j; \xx)$ into a 
single generating function $f(S; \xx)$.
We obtain the following result.
\proclaim{(3.10) Lemma} Let us fix $k$ and $d$. 
Then there exists a polynomial time 
algorithm, which, given polytopes $Q_1, \ldots, Q_n$ and functions 
$f(S \cap Q_j; \xx)$ computes $f(S; \xx)$ in the form
$$f(S; \xx)=\sum_{i \in I} \beta_i {\xx^{q_i} \over 
(1-\xx^{b_{i1}}) \cdots (1-\xx^{b_{is}})}$$
for $s \leq 2k$.
\endproclaim 
\demo{Proof}
We can write 
$$f(S; \xx)=\sum_{i=1}^n f(S \cap \inte Q_i; \xx).$$
On the other hand, 
$$S \cap \inte Q_i= (S \cap Q_i) \cap (\inte Q_i \cap {\Bbb Z}^d).$$
First, using Corollary 3.9 we compute $f(\inte Q_i \cap {\Bbb Z}^d; \xx)$,
and then using Theorem 3.6 we compute $f(S \cap \inte Q_i; \xx)$.
{\hfill \hfill \hfill} \qed
\enddemo

\head 4. Lattice Width and Small Gaps \endhead

In this section, we establish a simple geometric fact which plays a 
crucial role in the proof of Theorem 1.7.
We start with definitions.
\definition{(4.1) Definitions} 
Let $\Lambda \subset {\Bbb R}^d$ be a lattice (that is, a 
discrete additive subgroup of ${\Bbb R}^d$ of rank $d$) and
let $\Lambda^{\ast}  \subset {\Bbb R}^d$ be the dual (reciprocal) lattice,
that is, 
$$\Lambda^{\ast}=\Bigl\{c \in {\Bbb R}^d: \ \langle c, x \rangle \in {\Bbb Z}
\quad \text{for all} \quad x \in \Lambda \Bigr\},$$ 
where $\langle \cdot , \cdot \rangle$ is the standard scalar product in 
${\Bbb R}^d$.
For a convex body $B \subset {\Bbb R}^d$ (by which we mean a 
convex compact set) and a non-zero 
vector $c \in \Lambda^{\ast}$
let
$$\wi(B,c)=\max_{x \in B} \langle c, x \rangle -
\min_{x \in B} \langle c, x \rangle$$
be the {\it width} of $B$ in the direction of $c$. 
Let 
$$\wi(B)=\min_{c \in \Lambda^{\ast} \setminus \{0\}} \wi(B,c)$$ 
be the {\it lattice width} of $B$. As is well known, the minimum indeed 
exists.

It is known that there exists a constant $\omega(d)$ with the following 
property: if $B \cap \Lambda =\emptyset$ then $\wi(B) \leq \omega(d)$.
It is conjectured (and proved in many special cases) 
that $\omega(d)=O(d)$ while the best known value is
$\omega(d)=O(d \ln d)$ \cite{BL+99}. 
\enddefinition
We state some obvious properties of the width:
$$\split &\wi(B,c)=\wi(B+x, c) \quad \text{for any} \quad x \in {\Bbb R}^d 
\qquad \text{and} \\ &\wi(\alpha B, c)=\alpha \wi(B,c) \quad 
\text{for all} \quad \alpha \geq 0. \endsplit$$
Consequently, 
$$\split &\wi(B)=\wi(B+x) \quad \text{for any} \quad x \in {\Bbb R}^d 
\qquad \text{and} \\ &\wi(\alpha B)=\alpha \wi(B) \quad 
\text{for all} \quad \alpha \geq 0. \endsplit$$
\proclaim{(4.2) Lemma} Let $B \subset {\Bbb R}^d$ be a convex body, let 
$c \in {\Bbb R}^d$ be a non-zero vector and let 
$$\gamma_{\min}=\min_{x \in B} \langle c, x \rangle \quad \text{and}\quad 
\gamma_{\max}=\max_{x \in B} \langle c, x\rangle.$$
Let $\gamma_{\min} < \gamma_1 < \gamma_2 < \gamma_{\max}$ be numbers.
Then there exists a point $x_0 \in B$ and a number $0<\alpha <1$ 
such that for 
$$A=\alpha(B-x_0)+x_0=\alpha B + (1-\alpha) x_0$$ one has
$A \subset B$ and 
$$\min_{x \in A} \langle c, x \rangle =\gamma_1 \quad \text{and} \quad 
\max_{x \in A} \langle c, x\rangle =\gamma_2.$$
\endproclaim
\demo{Proof} Translating $B$, if necessary, we can assume that 
$\gamma_{\min}=0$. Dilating $B$, if necessary, we can assume that 
$\gamma_{\max}=1$. 
Then $0<\gamma_1/(1-\gamma_2+\gamma_1)<1$, and, therefore, we can 
choose $x_0 \in B$ such that $\langle c, x_0 \rangle=
\gamma_1/(1-\gamma_2+\gamma_1)$.
Let $\alpha=(\gamma_2-\gamma_1)$. Then, for 
$A=\alpha(B-x_0) +x_0=\alpha B +(1-\alpha) x_0$, we have
$$\min_{x \in A} \langle c, x \rangle=
{(1-\alpha) \gamma_1 \over 1-\gamma_2 +\gamma_1}=\gamma_1$$ 
and 
$$\max_{x \in A} \langle c, x \rangle =
\alpha +{(1-\alpha) \gamma_1 \over 1-\gamma_2 +\gamma_1}=\gamma_2.$$
Since $B$ is convex, we have $A \subset B$.
{\hfill \hfill \hfill} \qed
\enddemo
Now we can prove the main result of this section.
\proclaim{(4.3) Theorem} Let $B \subset {\Bbb R}^d$ be a convex body and let 
$\Lambda \subset {\Bbb R}^d$ be a lattice. 
Let $c \in \Lambda^{\ast}$ be a non-zero vector.
Let us consider the map:
$$\phi: B \cap \Lambda \longrightarrow {\Bbb Z}, \quad 
\quad \phi(x)=\langle c, x \rangle$$ 
and let $Y=\phi(B \cap \Lambda)$.
Hence $Y \subset {\Bbb Z}$ is a finite set.

Suppose that 
$$\wi(B,c) \leq 2 \wi(B).$$
Then for any $y_1, y_2 \in Y$ such that $y_2-y_1 >2 \omega(d)$ 
there exists 
a $y \in Y$ such that $y_1 <y< y_2$. 
\endproclaim
\demo{Proof} Suppose that such a point $y$ does not exist. 
Let us choose any $0<\epsilon<1/2$ and let 
$\gamma_1=y_1+\epsilon$ and $\gamma_2=y_2-\epsilon$.
By Lemma 4.2 there exists an $x_0 \in B$ and a number $\alpha >0$ such that
for $A=\alpha(B-x_0)+x_0$, $A \subset B$, we have 
$$\min_{x \in A} \langle c, x\rangle =\gamma_1 \quad \text{and} \quad 
\max_{x \in A} \langle c, x \rangle =\gamma_2.$$  
Then there is no integer in the interval $[\gamma_1, \gamma_2]$ 
which is a value of $\langle c, x \rangle$ for some 
$x \in B \cap \Lambda$. Hence
 $A \cap \Lambda=\emptyset$. Therefore, we must have 
$$\wi(A) \leq \omega(d).$$
On the other hand, since $A$ is a homothetic image of $B$ we have 
$$\wi(A)=\alpha \wi(B) \quad \text{and} \quad 
\wi(A,c)=\alpha \wi(B,c).$$
Therefore,
$$\gamma_2-\gamma_1=\wi(A,c) \leq 2 \wi(A) \leq 2 \omega(d).$$
Hence $y_2 -y_1 -2\epsilon \leq 2 \omega(d)$ for any $\epsilon>0$ and 
$y_2-y_1 \leq 2 \omega(d)$, which is a contradiction.
{\hfill \hfill \hfill} \qed   
\enddemo
In other words, the set $Y \subset {\Bbb Z}$ does not have ``gaps'' 
larger than $2 \omega(d)$. We will use the following corollary of Theorem 4.3.
\proclaim{(4.4) Corollary} 
Let $Y \subset {\Bbb Z}$ be the set of Theorem 4.3 and 
let $m= \lceil 2 \omega(d) \rceil$. For a positive integer $l$, let 
$Y+l=\bigl\{y+l:\ y \in Y \bigr\}$ denote the translation of $Y$ by $l$.
If $Y \ne \emptyset$ then the set 
$$Z=Y \setminus  \bigcup_{l=1}^m (Y+l)$$
consists of a single point. 
\endproclaim 
\demo{Proof} By Theorem 4.3, we have $Z=\{z\}$, where $z=\min\{y:\ y \in Y \}$.
{\hfill \hfill \hfill} \qed
\enddemo

\head 5. Projections and Partitions \endhead

In this section, we supply the remaining ingredient of the proof 
of Theorem 1.7. This ingredient, up to a change of the coordinates, is 
a weak form of a lemma of R. Kannan \cite{K92}. 

We describe it below.
Let $T: {\Bbb R}^d \longrightarrow {\Bbb R}^k$ be a linear transformation 
such that $T({\Bbb R}^d)={\Bbb R}^k$ and $T({\Bbb Z}^d) \subset {\Bbb Z}^k$.
Thus $k \leq d$ and the matrix of $T$ is integral with respect
to the standard bases 
of ${\Bbb R}^d$ and ${\Bbb R}^k$. Then $\ker T$ is a rational 
$(d-k)$-dimensional subspace of ${\Bbb R}^d$ (that is, a subspace 
spanned by integer vectors) and $\Lambda ={\Bbb Z}^d \cap (\ker T)$ is 
a lattice in $\ker T$. As is known (see, for example, 
Chapter 1 of \cite{C97}), a basis of $\Lambda$ can be 
extended to a basis of ${\Bbb Z}^d$ and hence any linear functional 
$\ell: \ker T \longrightarrow {\Bbb R}$ such that 
$\ell(\Lambda) \subset {\Bbb Z}$ can be represented in the form
$\ell(x)=\langle c, x \rangle$ for some $c \in {\Bbb Z}^d$. The 
representation, of course, is not unique as long as $\ker T \ne {\Bbb R}^d$. 
For $c \in (\ker T)^{\bot}$ (the orthogonal complement of $\ker T$), the 
corresponding linear functional is identically 0.

Let $P \subset {\Bbb R}^d$ be a rational polytope. For $y \in {\Bbb R}^k$ 
let us consider the fiber 
$$P_y=\Bigl\{x \in P: \quad T(x)=y \Bigr\}$$
of $x$. For $c \in {\Bbb Z}^d \setminus (\ker T)^{\bot}$ we define the 
width of $P_y$ in the direction of $c$ as 
$$\wi(P_y, c)=\max_{x \in P_y} \langle c, x \rangle -\min_{x \in P_y} 
\langle c, x \rangle$$ 
and we define the lattice width of $P_y$ as 
$$\wi(P_y)=\min_{c \in {\Bbb Z}^d \setminus (\ker T)^{\bot}} \wi(P_y, c).$$ 
We observe that the lattice width of $P_y$ so defined coincides with the 
width (as defined in Section 4), 
with respect to $\Lambda$ of a translation $P_y' \subset \ker T$.

We need the following result, which is a (rephrased) weaker version of 
Lemma 3.1 from \cite{K92}. It asserts, roughly, that one can
dissect the image $T(P)$ into polynomially many (in the input size of $P$ and
$T$) polyhedral 
pieces $Q_i$ and find for every piece $Q_i$ a 
lattice direction $w_i$ such that for all $y \in Q_i$ 
the lattice width of $P_y$ is almost attained at $w_i$.   
\proclaim{(5.1) Lemma} Let us fix $d$. Then there exists a polynomial time 
algorithm, which, for any rational polytope 
$P \subset {\Bbb R}^d$ and any linear transformation $T: {\Bbb R}^d
\longrightarrow {\Bbb R}^k$ such that 
$T({\Bbb R}^d)={\Bbb R}^k$ and
$T({\Bbb Z}^d)={\Bbb Z}^k$, constructs  
rational polytopes $Q_1, \ldots, Q_n \subset {\Bbb R}^k$ and vectors
$w_1, \ldots, w_n \in {\Bbb Z}^d \setminus (\ker T)^{\bot}$
such that 
\roster 
\item For each $i=1, \ldots, n$ and every $y \in Q_i$, 
$$\text{either} \quad \wi(P_y, w_i) \leq 1 \quad 
\text{or} \quad \wi(P_y, w_i) \leq 2 \wi(P_y);$$
\item The interiors $\inte Q_i$ are pairwise disjoint and 
$$\bigcup_{i=1}^n \inte Q_i = T(P).$$
\endroster 
\endproclaim
\demo{Proof} Let us construct a rational subspace $V \subset {\Bbb R}^d$ 
such that $V \cap (\ker T)=\{0\}$ and $(\ker T) + V ={\Bbb R}^d$. Then 
the restriction of $T$ onto $V$ is invertible and we can compute a 
matrix $L$ of 
the linear
transformation ${\Bbb R}^k \longrightarrow V$ which is the right 
inverse of $T$. 

Suppose that the polytope $P$ is defined by a system of 
linear inequalities 
$$P=\Bigl\{ x \in {\Bbb R}^d: \quad Ax \leq b \Bigr\},$$
where $A$ is an $n \times d$ integer matrix and $b$ is an integer $n$-vector.
Then the translation $P_y' \subset \ker T$ of $P_y$ is defined by 
the system of linear inequalities
$$P_y'=\Bigl\{x \in \ker T:\quad  Ax \leq b-ALy \Bigr\}.$$
As $y$ ranges over $Q=T(P)$, vector $b'=b-ALy$ ranges over 
the rational polytope $Q'=b-AL(Q)$ with 
$\dim Q' \leq k$. Since $\wi(P_y, c)=\wi(P_y', c)$ for all $y \in Q$ and 
all $c$ and $\wi(P_y)=\wi(P_y')$, the result follows by Part 3 of Lemma 
3.1 of \cite{K92}. 
{\hfill \hfill \hfill} \qed
\enddemo

\head 6. Proofs \endhead

Now we are ready to prove Theorem 1.7.
\demo{Proof of Theorem 1.7} Without loss of generality, we assume 
that $T({\Bbb R}^d)={\Bbb R}^k$. Indeed, if $\im(T) \ne {\Bbb R}^k$,
we consider the restriction $T:{\Bbb R}^d \longrightarrow \im(T)$.
After a change of the coordinates, the lattice 
$\Lambda={\Bbb Z}^k \cap \im(T)$ is identified with the standard 
integer lattice. 

The proof is by induction on
$\dim (\ker T)=d-k$.

Suppose that $k=d$, so $\dim (\ker T)=0$ and 
$T: {\Bbb Z}^d \longrightarrow {\Bbb Z}^k={\Bbb Z}^d$ is 
an embedding. Let $e_1, \ldots, e_d$ be the standard basis of 
${\Bbb Z}^d$ and let $t_i=T(e_i)$. Then 
$f(S; \xx)$ is obtained from $f(P \cap {\Bbb Z}^d; \yy)$ by 
the monomial substitution $y_i=\xx^{t_i}$ and we use Theorems 2.6
and 3.1 to complete the proof.

Suppose that $d>k$, so $\dim (\ker T)>0$. Let  
$Q_1, \ldots, Q_n \subset {\Bbb R}^k$ be the polytopes constructed in 
Lemma 5.1.
It suffices to compute
the functions $f(S \cap Q_i; \xx)$ for $i=1, \ldots, n$ and 
then, using Lemma 3.10, we can patch them together and obtain 
$f(S; \xx)$.

Let us consider a particular polytope $Q=Q_i$ 
and the corresponding intersection
$S \cap Q$. Let $w=w_i$, $w \in {\Bbb Z^d} \setminus (\ker T)^{\bot}$ 
be a vector whose existence is claimed by Lemma 5.1.
Let us consider the linear transformation 
$$\hat{T}: {\Bbb R}^d \longrightarrow {\Bbb R}^{k+1}={\Bbb R}^k \oplus 
{\Bbb R}, \quad
\hat{T}(x)=\bigl(T(x),\ \langle w, x \rangle\bigr)$$
and the projection 
$$pr: {\Bbb R}^{k+1} \longrightarrow {\Bbb R}^k, \quad 
pr(\xi_1, \ldots, \xi_{k+1})=(\xi_1, \ldots, \xi_k).$$
Finally, let
$P'=\bigl\{x \in P:\ T(x) \in Q \bigr\}$ and
 $\hat{S}=\hat{T}(P' \cap {\Bbb Z}^d) \subset {\Bbb R}^{k+1}$. 

Clearly, $S\cap Q=pr(\hat{S})$ and $\dim(\ker \hat{T})=d-k-1$, so we can 
apply the induction hypothesis to $\hat{T}$ and compute 
$f(\hat{S}; \zz)$, where $\zz=(\xx, x_{k+1})$, $x_{k+1} \in {\Bbb C}$. 
Our goal is to compute $f(S \cap Q; \xx)$ 
from $f(\hat{S}; \zz)$. To do that,
we construct a subset $Z \subset \hat{S}$ such that the projection
$pr: Z \longrightarrow S \cap Q$ is one-to-one, and then we obtain 
$f(S \cap Q; \xx)$ from 
$f(Z; \zz)$ by substituting $x_{k+1}=1$.

For a positive integer $l$, let $\hat{S}+l$ denote the translation 
of $\hat{S}$ by $l$ along the last coordinate,
$$\hat{S}+l=\bigl\{(\xi_1, \ldots, \xi_k, \xi_{k+1}+l): \quad
(\xi_1, \ldots, \xi_k) \in \hat{S} \bigr\}.$$
Clearly,
$$f({\hat S}+l; \zz)=x_{k+1}^l f(\hat{S}; \zz).$$
Let $m=\lceil 2 \omega(d-k) \rceil$ (see Section 4) and let us define
$$Z=\hat{S} \setminus \bigcup_{l=1}^m \bigl( \hat{S}+l \bigr).$$ 
Using Corollary 3.7, we compute $f(Z; \zz)$. 

Now we claim that the 
projection $pr: Z \longrightarrow S \cap Q$ is one-to-one.
Let us 
consider the projection $pr: {\hat S} \longrightarrow S \cap Q$. 
For a $y \in S$ let us consider 
the preimage $\hat{S}_y \subset \hat{S}$ of $y$.
We observe that 
$$\hat{S}_y=
\Bigl\{\bigl(y,\ \langle w, x \rangle\bigr): 
\quad x \in P_y \cap {\Bbb Z}^d\Bigr\},$$
that is, $\hat{S}_y$ consists of all pairs 
$\bigl(y,\ \langle w, x \rangle \bigr)$, where $x$ is an integer point 
from the fiber $P_y$ of $P$ over $y$:
$$P_y=\Bigl\{x \in P: \quad T(x)=y \Bigr\}.$$ 
By Lemma 5.1, we have either 
$\wi(P_y, w) \leq 1$ or $\wi(P_y, w) \leq 2 \wi(P_y)$.
If $\wi(P_y, w) \leq 2 \wi(P_y)$, then, by Corollary 4.4,
the set 
$$Z_y=\hat{S}_y \setminus \bigcup_{l=1}^m \bigl(\hat{S}_y + l)$$ 
consists of a single point, that is, the point of $\hat{S}_y$ with the 
smallest last coordinate.   
If $\wi(P_y, w) \leq 1$ then $\hat{S}_y$ consists of a single point 
and so $Z_y$ consists of a single point as well.
Thus, in any case, for any $y \in S \cap Q$ the preimage 
$Z_y$ of the projection $pr: Z \longrightarrow S \cap Q$ consists of a single point,
so $pr: Z \longrightarrow S \cap Q$ is indeed one-to-one. Hence, using 
Theorem 2.6, we compute $f(S \cap Q; \xx)$ by specializing $f(Z; \zz)$ at
$x_{k+1}=1$ (where $\zz=(\xx, x_{k+1})$).
{\hfill \hfill \hfill} \qed  
\enddemo
We deduce Theorem 1.5 from Theorem 1.7.
\demo{Proof of Theorem 1.5}
Let us define a linear transformation
 $T: {\Bbb R}^d \longrightarrow {\Bbb R}$ by 
$$T(\xi_1, \ldots, \xi_d)=a_1 \xi_1 + \ldots + a_d \xi_d.$$ Thus 
$S=T({\Bbb Z}^d_+)$ is the semigroup generated by $a_1, \ldots, a_d$. 
It remains to notice that there are some explicit bounds for the largest 
positive integer not in $S$, so one can replace 
the non-negative orthant ${\Bbb Z}^d_+$ by a rational polytope to get 
the initial interval of $S$.
 For example, in \cite{EG72} it is 
shown that if $t \geq \max\{a_1, \ldots, a_d\}$ 
then all numbers greater than 
or equal to $2t^2/d$ are in $S$. Let $n=\lceil 2t^2/d \rceil$ and let 
$$P=\Bigl\{(\xi_1, \ldots, \xi_d):\quad \sum_{i=1}^d \xi_i a_i \leq n-1
\quad \text{and} \quad
\xi_i \geq 0 \quad 
\text{for} \quad i=1, \ldots d \Bigr\}$$ be the simplex
 in ${\Bbb R}^d$. Then we can represent 
$S$ as a disjoint union of $T(P \cap {\Bbb Z}^d)$ and the integer points 
in the ray $[n, +\infty)$
Since the generating function of the set of 
integer points in the ray $[n, +\infty)$ is just $x^{n+1}/(1-x)$,
applying Theorem 1.7 we complete the proof.
{\hfill \hfill \hfill} \qed
\enddemo

\head 7. Further Examples: Hilbert Bases, Test Sets and 
Hilbert Series \endhead

As another application of Theorem 1.7, let us show that certain 
{\it Hilbert bases} are enumerated by short rational functions.

Let $u_1, \ldots, u_d \subset {\Bbb Z}^d$ be linearly independent vectors and
let 
$$\Pi=\Bigl\{ \sum_{i=1}^d \alpha_i u_i: \quad 0 \leq \alpha_i \leq 1 
\quad \text{for} \quad i=1,\ldots,d
\Bigr\}$$ 
be the parallelepiped spanned by $u_1, \ldots, u_d$, 
and let $K$ be the convex cone spanned by $u_1, \ldots, u_d$:
$$K=\Bigl\{\sum_{i=1}^d \alpha_i u_i: \quad \alpha_i \geq 0 \quad 
\text{for} \quad i=1, \ldots, d \Bigr\}.$$  
 We say that 
a point $v \in \Pi \cap {\Bbb Z}^d$ is {\it indecomposable} provided 
$v$ cannot be written in the form $v=v_1+v_2$, where $v_1$ and $v_2$ 
are non-zero integer points from $\Pi$. 
The set $S$ of all non-decomposable integer vectors in $\Pi$
is called the (minimal) {\it Hilbert basis} of the semigroup 
$K \cap {\Bbb Z}^d$, since every integer vector in $K$ can be written 
as a non-negative integer combination of points from $S$, see Section 16.4 of
\cite{Sc86}.
Let us show that as long as the dimension $d$ is fixed, the set $S$ 
has a short rational generating function.
\proclaim{(7.1) Theorem} Let us fix $d$. Then
there exists a number $s=s(d)$ and
 a polynomial time algorithm, which, given linearly 
independent vectors $u_1, \ldots, u_d \in {\Bbb Z}^d$ computes 
the generating function $f(S; \xx)$ of the (minimal) Hilbert basis $S$
of the semigroup of integer points in the cone spanned by 
$u_1, \ldots, u_d$ in the form: 
$$f(S; x)=\sum_{i \in I} \alpha_i {x^{p_i} \over (1-x^{b_{i1}}) \cdots 
(1-x^{b_{is}})},$$
where $I$ is a set of indices, $\alpha_i$ are rational numbers, 
$p_i$, $b_{ij} \in {\Bbb Z}^d$ and $b_{ij} \ne 0$ for all $i,j$.
\endproclaim
\demo{Proof} Let us construct a rational polyhedron $Q \subset \Pi$ which 
contains all integer points in $\Pi$ except 0. This can be done, for 
example, as follows: we construct
 vectors $l_1, \ldots, l_d \in {\Bbb Z}^d$ such that 
$\langle l_i, u_j \rangle =0$ for $i \ne j$ and $\langle l_i, u_i \rangle >0$,
let $l=l_1 + \ldots + l_d$ and intersect $\Pi$ with the halfspace 
$\langle l, x \rangle \geq 1$. 

Let $P=Q \times Q \subset {\Bbb R}^d \oplus {\Bbb R}^d={\Bbb R}^{2d}$ and 
let $T: P \longrightarrow {\Bbb R}^d$ be the transformation, 
$T(x,y)=x+y$. Let $S_1=T(P \cap {\Bbb Z}^{2d})$ and let 
$S_2=Q \cap {\Bbb Z}^d$. Then the minimal Hilbert base $S$ can be 
written as $S=S_2 \setminus S_1$. The proof now follows from Theorem 1.7 and 
Corollary 3.7.
{\hfill \hfill \hfill} \qed 
\enddemo
Yet another interesting class of sets having short rational generating 
functions is that of ``test sets'' with respect to a
given integer matrix.
\subhead (7.2) Test Sets \endsubhead
Let us choose a $n \times d$ integer matrix $A$ such that for any 
$b \in {\Bbb R}^n$, the polyhedron
$$P_b=\Bigl\{x \in {\Bbb R}^d: \quad Ax \leq b \Bigr\},$$
is bounded.
A point $a \in {\Bbb Z}^d$ is called a {\it neighbor} of $0$ with 
respect to $A$ provided there is a polytope $P_b$ containing $0$ and $a$
and not containing any other integer point in its interior.
The set $S(A)$ of all neighbors of the origin is often called a {\it test
set}. 
Test sets $S(A)$ play an important role in parametric integer programming 
\cite{S97}. The set $S(A)$ is finite, and it has some interesting 
(for $d \geq 2$) and not quite understood 
(for $d \geq 3$) structure. One can show 
that for any fixed $d$ and $n$, given $A$, 
the generating function $f(S; \xx)$ for 
$S=S(A)$ can be computed in polynomial time as a short rational
function. The proof follows from Theorem 7.1 and Corollary 3.7 in a similar 
way as above, since $S$ can be expressed as a Boolean combination of 
projections of sets of integer points in some rational polytopes.

We note that other types of test sets studied in the literature, 
such as Schrijver's universal test 
set and Graver's test set (see \cite{T95} and \cite{St96}), 
also admit a short rational generating function.
\bigskip
Finally, we describe one related problem of computational commutative 
algebra.
\subhead (7.3) Hilbert series of rings generated by monomials \endsubhead
Let us consider integer vectors $a_1, \ldots, a_d \in {\Bbb Z}^k_+$ 
with non-negative coordinates
and let $S$ be the semigroup generated by $a_1, \ldots, a_d$:
$$S=\Bigl\{ \sum_{i=1}^d \mu_i a_i: \quad \mu_i \in {\Bbb Z}_+ \Bigr\}.$$
Thus $S$ can be represented as the image $T({\Bbb Z}^d_+)$ under the 
linear transformation 
$$T: {\Bbb R}^d \longrightarrow {\Bbb R}^k, \quad 
T(\xi_1, \ldots, \xi_d)=\xi_1 a_1 + \ldots + \xi_d a_d.$$
The generating function $f(S;\xx)$ can be interpreted as the Hilbert 
series of the ${\Bbb Z}^k$-graded ring $R={\Bbb C}[\xx^{a_1}, \ldots, 
\xx^{a_d}]$, cf. \cite{BS98} and Chapter 10 of \cite{St96}.
The set $S$ is infinite and Theorem 1.7 is not directly applicable
(although it allows us to claim the intersection of $S$ with any given 
polytopal region has a short rational generating function).
However, one can still compute the whole function
$f(S;\xx)$ in polynomial time as a 
short rational function 
provided the dimension $k$ and the number $d$ of generators are fixed. 
We also note that by applying a 
monomial specialization of
$f(S; \xx)$ we can obtain the Hilbert series of $R$ under a coarser grading. 

We sketch an algorithm for computing $f(S; \xx)$ below.

Without loss of generality we assume that $a_i \ne 0$ for 
$i=1, \ldots, d$.
Let us consider the product 
$$g(S; \xx)=f(S; \xx) (1-\xx^{a_1}) \cdots (1-\xx^{a_d}).$$ 
It is not hard 
to prove that $g(S; \xx)$ is, in fact, a polynomial in $\xx$.
This follows, for example, from the interpretation of $f(S;\xx)$ 
as a Hilbert series, cf. Section I.9 of \cite{E95}.
 
We need to compute a bound $L$ with the property that that 
if the coefficient of $\xx^m$, $m=(\mu_1, \ldots, \mu_k)$,
in $g(S; \xx)$ is non-zero then 
$\mu_1 + \ldots + \mu_k \leq L$.
Suppose for a moment that we can find 
such an $L$. Let us consider the integer cube 
$$C \subset {\Bbb Z}_+^d, \quad C=\Bigl\{ (\xi_1, \ldots, \xi_d): \quad 
0 \leq \xi_i \leq L \quad \text{for} \quad  i=1, \ldots,d  \Bigr\}$$
and the integer simplex
$$\Delta \subset {\Bbb Z}_+^k: \quad 
\Delta=\bigl\{ (\mu_1, \ldots, \mu_k): \quad 
\mu_1 + \ldots +\mu_k \leq L \bigr\}.$$
Let $S'=T(C)$, so $\Delta \cap S \subset S'$. 
Applying Theorem 1.7, we compute $f(S'; \xx)$ as a short rational 
function. Let 
$$g(S'; \xx)=f(S'; \xx) (1-\xx^{a_1}) \cdots (1-\xx^{a_d}).$$
We note that 
$$g(S; \xx)=g(S'; \xx) \star f(\Delta; \xx)$$ 
and use Theorem 3.1 and Lemma 3.4 to compute 
the Hadamard product $g(S; \xx)$ as a short 
rational function. Finally,
we let 
$$f(S; \xx)=g(S; \xx) \prod_{i=1}^k {1 \over 1-\xx^{a_i}}.$$ 

It remains, therefore, to compute the bound 
$L$ on the total degree of a monomial $\xx^m$ which may appear 
with a non-zero coefficient in the expansion of $g(S; \xx)$.

Let us consider the rational cone $K \subset {\Bbb R}^d \oplus 
{\Bbb R}^d$,
$$K=\Bigl\{(x,y): \quad x, y \in {\Bbb R}^d_+ \quad \text{and} \quad 
T(x)=T(y) \Bigr\}.$$
The lattice semigroup $K \cap {\Bbb Z}^{2d}$ is finitely generated
and using some standard techniques (see Chapter 17 of 
\cite{Sc86} and Chapter 4 of \cite{St96}) one can 
compute in polynomial time an upper bound $M$ on the coordinates of 
generators $(x_i, y_i)$ of 
$K \cap {\Bbb Z}^{2d}$. Let $A$ be the sum of the coordinates 
of $a_1, \ldots, a_d$.
We claim that $L=A(M+1)$ is the desired upper bound.
  
Indeed, for every generator $(x_i, y_i)$ with $x_i \ne y_i$, let
$z_i=x_i-y_i$ or $z_i=y_i-x_i$, whichever is lexicographically
positive. Thus each coordinate of $z_i$ does not exceed $M$.  
Let 
$$Z={\Bbb Z}^d_+ \setminus \bigcup_i \bigl({\Bbb Z}^d_+ + z_i\bigr).$$
One can observe that the restriction $T: Z \longrightarrow S$ is 
one-to-one. In fact, for every $x \in S$ the vector 
$z \in Z$ such that $T(z)=x$ is the lexicographic minimum among 
all $y \in {\Bbb Z}^d_+$ such that $T(y)=x$.

For $I \subset \{1, \ldots, d\}$ let ${\Bbb Z}^I_+ \subset 
{\Bbb Z}^d_+$ be the coordinate semigroup consisting of the points 
$(\xi_1, \ldots, \xi_d)$ such that $\xi_i=0$ for $i \notin I$. 
As is proved in \cite{Kh95}, the set $Z$ can be represented as 
a finite disjoint union of sets $Z_j$ of the type 
$v_j + {\Bbb Z}^{I_j}_+$ so that the coordinates of $v_j$ do not exceed
$M$. Let $S_j=T(Z_j)$. Then $S$ is the disjoint union of $S_j$ and 
$$f(S_j; \xx)=\xx^{T(v_j)} \prod_{i \in I_j} {1 \over 1-\xx^{a_i}}.$$
The sum of the coordinates of $T(v_j)$ does not exceed $MA$. Therefore,
if $\xx^m$, $m=(\mu_1, \ldots, \mu_k)$ appears with a non-zero coefficient 
in the product 
$$f(S_j; \xx) (1-\xx^{a_1}) \cdots (1-\xx^{a_k}),$$ 
we must have $\mu_1 + \ldots +\mu_d \leq  MA+A=L$, which completes 
the proof.

\head References \endhead
\refstyle{A}


\widestnumber\key{AAAAA}
\Refs\nofrills{}

\ref\key{B94}
\by A. Barvinok
\paper  A polynomial time algorithm for counting integral points in polyhedra when the dimension is fixed
\jour Math. Oper. Res. 
\vol 19
\yr 1994
\pages 769--779
\endref

\ref\key{B02}
\by A. Barvinok
\book A Course in Convexity
\bookinfo Graduate Studies in Mathematics
\vol 54
\publ Amer. Math. Soc.
\publaddr Providence, RI
\yr 2002
\endref

\ref \key{BP99}
\by A. Barvinok and J.E. Pommersheim
\paper An algorithmic theory of lattice points in polyhedra
\inbook New Perspectives in Algebraic Combinatorics
 (Berkeley, CA, 1996--97)
\pages  91--147 
\bookinfo Math. Sci. Res. Inst. Publ.
\vol 38
\publ Cambridge Univ. Press
\publaddr Cambridge
\yr 1999
\endref

\ref \key{BS98}
\by D. Bayer and B. Sturmfels
\paper Cellular resolutions of monomial modules
\jour J. Reine Angew. Math.
\vol 502 
\yr 1998
\pages 123--140
\endref

\ref \key{BL+99}
\by W. Banaszczyk, A.E. Litvak, A. Pajor and S.J. Szarek
\paper The flatness theorem for nonsymmetric convex bodies via the local
theory of Banach spaces
\jour Math. Oper. Res.
\vol 24 
\yr 1999
\pages 728--750
\endref 

\ref \key{C97}
\by J.W.S Cassels
\book An Introduction to the Geometry of Numbers. 
Corrected reprint of the 1971 edition
\bookinfo Classics in Mathematics
\publ Springer-Verlag
\publaddr Berlin
\yr 1997
\endref

\ref\key{D96}
\by G. Denham
\paper The Hilbert series of a certain module
\paperinfo manuscript
\yr 1996
\endref

\ref\key{E95}
\by D. Eisenbud
\book Commutative Algebra with a View Toward Algebraic Geometry
\bookinfo Graduate Texts in Mathematics
\vol 150
\publ Springer-Verlag
\publaddr New York
\yr 1995
\endref

\ref \key{EG72}
\by P. Erd\"os and R.L. Graham 
\paper On a linear diophantine problem of Frobenius
\jour Acta Arith.
\vol 21 
\yr 1972
\pages 399--408
\endref

\ref\key{GLS93}
\by M. Gr\"otschel, L. Lov\'asz and A. Schrijver
\book Geometric Algorithms and Combinatorial Optimization. Second edition
\bookinfo Algorithms and Combinatorics
\vol 2
\publ Springer-Verlag
\publaddr Berlin
\yr 1993
\endref

\ref \key{K92}
\by R. Kannan
\paper Lattice translates of a polytope and the Frobenius problem
\jour Combinatorica 
\vol 12 
\yr 1992
\pages 161--177
\endref

\ref \key{Kh95}
\by A.G. Khovanskii
\paper Sums of finite sets, orbits of commutative semigroups and Hilbert
functions (Russian)
\jour Funktsional. Anal. i Prilozhen. 
\vol 29 
\yr 1995
\pages 36--50
\transl translated in {\it Funct. Anal. Appl.}, {\bf  29}(1995), 
 102--112.
\endref


\ref \key{KLS90}
\by R. Kannan, L. Lov\'asz, and H. Scarf
\paper The shapes of polyhedra
\jour Math. Oper. Res. 
\vol 15 
\yr 1990
\pages 364--380
\endref

\ref\key{P94}
\by C.H. Papadimitriou
\book Computational Complexity
\publ Addison-Wesley 
\publaddr Reading, MA 
\yr 1994
\endref

\ref \key{S97}
\by H. Scarf
\paper Test sets for integer programs
\jour Math. Programming, Ser. B 
\vol 79 
\yr 1997
\pages 355--368
\endref

\ref \key{Sc86}
\by A. Schrijver
\book Theory of Linear and Integer Programming
\publ  Wiley-Interscience 
\publaddr Chichester
\yr 1986
\endref

\ref \key{St96}
\by B. Sturmfels
\book Gr\"obner Bases and Convex Polytopes
\bookinfo University Lecture Series
\vol 8
\publ Amer. Math. Soc.
\publaddr Providence, RI
\yr 1996
\endref

\ref \key{SW86}
\by L.A. Sz\'ekely and N.C. Wormald
\paper Generating functions for the Frobenius problem with $2$ and $3$ 
generators
\jour Math. Chronicle 
\vol 15 
\yr 1986
\pages 49--57
\endref

\ref \key{T95}
\by R. Thomas
\paper A geometric Buchberger algorithm for integer programming
\jour Math. Oper. Res.
\vol 20
\yr 1995 
\pages 864--884
\endref



\endRefs
\enddocument
\end